\documentclass[a4paper,12pt]{article}
\usepackage{amsthm,amsmath}
\usepackage{amssymb}
\usepackage[all]{xypic}
\usepackage{enumerate}
\usepackage{a4wide}
\usepackage{hyperref}

\swapnumbers
\newtheorem{theorem}{Theorem}[section]
\newtheorem{proposition}[theorem]{Proposition}
\newtheorem{corollary}[theorem]{Corollary}
\newtheorem{lemma}[theorem]{Lemma}
\newtheorem*{theorem*}{Theorem}
\newtheorem*{proposition*}{Proposition}
\newtheorem*{corollary*}{Corollary}
\newtheorem*{lemma*}{Lemma}
\theoremstyle{definition}
\newtheorem{definition}[theorem]{Definition}

\newtheorem{example}[theorem]{Example}
\newtheorem{remark}[theorem]{Remark}
\newtheorem*{remark*}{Remark}

\newtheorem*{definition*}{Definition}

\newcommand{\coring}[1]{\mathfrak{#1}}
\newcommand{\tensor}[1]{\otimes_{#1}}
\newcommand{\tensfun}[1]{\underset{{#1}}{\otimes}}

\newcommand{\rcomod}[1]{\mathcal{M}^{#1}}
\newcommand{\rmod}[1]{\mathcal{M}_{#1}}
\newcommand{\lmod}[1]{{}_{#1}\mathcal{M}}

\newcommand{\cotensor}[1]{\square_{#1}}
\newcommand{\lcomod}[1]{{}^{#1}\mathcal{M}}

\renewcommand{\hom}[3]{\mathrm{Hom}_{#1}(#2,#3)}

\newcommand{\rend}[2]{\mathrm{End}({#2}_{#1})}

\newcommand{\lend}[2]{\mathrm{End}({}_{#1}#2)}

\newcommand{\cohom}[3]{\mathrm{h}_{#1}(#2,#3)}

\newcommand{\rcomatrix}[2]{#2^* \tensor{#1} #2}

\begin{document}
\title{Comatrix corings: Galois corings, Descent Theory, and a Structure Theorem for Cosemisimple corings }
\author{L. El Kaoutit {\normalsize and} J. G\'omez-Torrecillas \\
\normalsize Departamento de \'{A}lgebra \\ \normalsize Universidad de
Granada\\ \normalsize E18071 Granada, Spain \\ \normalsize e-mail:
\textsf{kaoutit@fedro.ugr.es} \\ \normalsize e-mail:
\textsf{torrecil@ugr.es} }

\date{}

\maketitle

\section*{Introduction}

Corings and their comodules provide an appropriate formalism to
unify notions and results coming from different subfields of the
theory of associative algebras. One of the most striking examples
is the following.  Let $\psi : T \rightarrow A$ be a ring
extension. M. Cipolla \cite{Cipolla:1976} extended Grothendieck's
theory of the faithfully flat descent from the commutative case to
the non commutative one. His main result says, in the restatement
given in \cite[Theorem 3.8]{Nuss:1997}, that if ${}_TA$ is
faithfully flat, then the tensor product functor $- \tensor{T} A :
\rmod{T} \rightarrow Desc_{\psi}$ establishes an equivalence
between the category $\rmod{T}$ of all right $T$--modules and the
category $Desc_{\psi}$ of descent data. Assume now that $A$ is a
right comodule algebra over a Hopf algebra $H$ and that $T$ is the
subring of coinvariant elements of $A$ (see
\cite{Schneider:1990a}). H. J. Schneider proved that the functor
$- \tensor{A} : \rmod{T} \rightarrow \mathcal{M}_A^H$, where
$\mathcal{M}_A^H$ is the category of right Hopf $A$--modules, is
an equivalence if and only if ${}_TA$ is faithfully flat and the
canonical map $can : A \tensor{T} A \rightarrow A \tensor{} H$ is
bijective \cite[Theorem 1]{Schneider:1990a}. Both theorems have a
similar flavor, and in fact they are particular cases of a recent
result on corings having a grouplike element due to T.
Brzezi\'{n}ski \cite[Theorem 5.6]{Brzezinski:2000unp}. The
enlightening fact here is that the categories $Desc_{\psi}$ and
$\mathcal{M}_A^H$ are categories of comodules over suitable
corings which become isomorphic to Sweedler's canonical corings of
the form $A \tensor{T} A$ (see \cite{Sweedler:1975}).  From a
categorical point of view, the $A$--corings characterized in
\cite[Theorem 5.6]{Brzezinski:2000unp} are precisely those for
which $A$ has a structure of right comodule such that $A$ becomes
a progenerator for the category of right comodules (see
\cite[Theorem 2.4]{Abuhlail:2002unp} and \cite[Theorem
3.5]{Caenepeel/Vercruysse/Wang:2002unp}, for the finitely and
projective generated case). We think that the theory of corings
should be developed by its own right, so, in the light of the
aforementioned results, a basic question arising here is to
characterize the corings whose category of comodules has a
finitely generated projective generator. In order to extricate the
structure of these corings, we introduce the notion of comatrix
coring, which allows also to give a full generalization of
Brzezi\'nski's result which works for corings without grouplike
elements. Comatrix corings allow as well to give a complete
description of all cosemisimple corings.

We have organized our exposition as follows. After a brief
introduction (Section \ref{preliminaries}) to the basic notions of
corings and comodules, we expound in Section \ref{comatrixcan} how
to construct a comatrix $A$--coring $\Sigma^* \tensor{B} \Sigma$
from a given bimodule ${}_B\Sigma_A$ such that $\Sigma_A$ is
finitely generated and projective. In the case that $\Sigma$ is a
right comodule over an $A$--coring $\coring{C}$, we define a
homomorphism of $A$--corings $can : \Sigma^* \tensor{B} \Sigma
\rightarrow \coring{C}$. This canonical morphism generalizes the
homonymous map introduced in \cite{Schneider:1990a}.

Section \ref{gl-dethy} contains several characterizations of those
corings $\coring{C}$ having a finitely generated and projective
generator. We prove in particular that they are comatrix corings
$\Sigma^* \tensor{T} \Sigma$ such that ${}_T\Sigma$ is faithfully
flat. We include some consequences, among them, we deduce
\cite[Theorem 5.6]{Brzezinski:2000unp} and a Descent Theorem for
ring extensions of the form $T \rightarrow \rend{A}{\Sigma}$ which
generalizes \cite[Teorema]{Cipolla:1976} and \cite[Theorem
3.8]{Nuss:1997}.

In Section \ref{cosemisimple} we offer a structure theorem for
cosemisimple corings. They are described in a unique way as direct
sums of comatrix corings of the form $\Sigma^* \tensor{D} \Sigma$,
where $D$ is a division subring of $\rend{A}{\Sigma}$.

Section \ref{coend} is devoted to show that comatrix corings can
be alternatively introduced as coendomorphism corings.

\section{The basic notions}\label{preliminaries}

We use the following conventions. For an object $C$ in a category
the identity morphism $C \rightarrow C$ is denoted by $C$.  We
work over fixed commutative ring $K$, and all our additive
categories are assumed to be $K$--linear. For instance, all rings
in this paper are unitary $K$--algebras, and all bimodules are
assumed to centralize the elements of $K$. Attached to every
object $C$ of an additive category $\mathcal{A}$ we have its
endomorphism ring $\mathrm{End}_{\mathcal{A}}(C)$, whose
multiplication is given the composition of the category. As usual,
some special conventions will be understood for the case of
endomorphism rings of modules. Thus, if $M_R$ is a right module
over a ring $R$, then its endomorphism ring in the category
$\rmod{R}$ of all right $R$--modules will be denoted by
$\rend{R}{M}$, while if ${}_RN$ is a left $R$--module, then its
endomorphism ring, denoted by $\lend{R}{N}$, is, by definition,
the opposite of the endomorphism ring of $N$ in the category
$\lmod{R}$ of all left modules over $R$.

Throughout this paper, $A, A', \dots, B, \dots$ denote associative
and unitary algebras over a commutative ring $K$. The tensor
product over $A$ is denoted by $\tensor{A}$. We shall sometimes
replace $\tensor{K}$ by $\tensor{}$.

We recall from \cite{Sweedler:1975} the notion of a coring. An
$A$--\emph{coring} is a three-tuple ($\coring{C},
\Delta,\epsilon$) consisting of an $A$-bimodule $\coring{C}$ and
two $A$--bimodule maps
\[
\Delta : \coring{C} \longrightarrow \coring{C} \tensor{A}
\coring{C}, \qquad \epsilon : \coring{C} \longrightarrow A
\]
such that $(\coring{C} \tensor{A} \Delta) \circ \Delta = (\Delta
\tensor{A} \coring{C}) \circ \Delta$ and $(\epsilon \tensor{A}
\coring{C}) \circ \Delta = ( \coring{C} \tensor{A} \epsilon) \circ
\Delta = \coring{C}$. A \emph{right} $\coring{C}$--\emph{comodule}
is a pair $(M,\rho_M)$ consisting of right $A$--module $M$ and an
$A$--linear map $\rho_M: M \rightarrow M \tensor{A} \coring{C}$
satisfying $(M \tensor{A} \Delta) \circ \rho_M = (\rho_M
\tensor{A} \coring{C}) \circ \rho_M$, $(M \tensor{A} \epsilon)
\circ \rho_M = M$; such $M$ will be denoted by $M_{\coring{C}}$. A
\emph{morphism} of right $\coring{C}$--comodules $(M,\rho_M)$ and
$(N,\rho_N)$ is a right $A$--linear map $f: M \rightarrow N$ such
that $(f \tensor{A} \coring{C}) \circ \rho_M = \rho_N \circ f$;
the $K$--module of all such morphisms will be denoted by
$\hom{\coring{C}}{M}{N}$. The right $\coring{C}$--comodules
together with their morphisms form the additive category
$\rcomod{\coring{C}}$. Coproducts and cokernels in
$\rcomod{\coring{C}}$ do exist and can be already computed in
$\rmod{A}$. Therefore, $\rcomod{\coring{C}}$ has arbitrary
inductive limits. If ${}_A\coring{C}$ is flat, then
$\rcomod{\coring{C}}$ an abelian category. The converse is not
true, as \cite[Example 1.1]{ElKaoutit/Gomez/Lobillo:2001pr} shows.

 Let $\rho_M : M
\rightarrow M \tensor{A} \coring{C}$ be a comodule structure over
an $A'-A$--bimodule $M$, and assume that $\rho_M$ is $A'$--linear.
For any right $A'$--module $X$, the right $A$--linear map $X
\tensor{A'} \rho_{M} : X \tensor{A'} M \rightarrow X \tensor{A'} M
\tensor{A} \coring{C}$ makes $X \tensor{A'} M$ a right
$\coring{C}$--comodule. This leads to an additive functor $-
\tensor{A'} M : \rmod{A'} \rightarrow \rcomod{\coring{C}}$. The
classical adjointness isomorphism $\hom{A}{Y \tensor{B}M}{X} \cong
\hom{B}{Y}{\hom{A}{M}{X}}$ induces, by restriction, a natural
isomorphism $\hom{\coring{C}}{Y \tensor{B}M}{X} \cong
\hom{B}{Y}{\hom{\coring{C}}{M}{X}}$, for $Y \in \rmod{A'}, X \in
\rcomod{\coring{C}}$. Therefore, $\hom{\coring{C}}{M}{-} :
\rcomod{\coring{C}} \rightarrow \rmod{A'}$ is right adjoint to $-
\tensor{A'} M : \rmod{A'} \rightarrow \rcomod{\coring{C}}$. On the
other hand, the functor $- \tensor{A} \coring{C}$ is right adjoint
to the forgetful functor $U : \rmod{A} \rightarrow
\rcomod{\coring{C}}$ (see \cite[Proposition 3.1]{Guzman:1989},
\cite[Lemma 3.1]{Brzezinski:2000unp}).

Now assume that the $A'-A$--bimodule $M$ is also a left
$\coring{C}'$--comodule with structure map $\lambda_M : M
\rightarrow \coring{C}' \tensor{A} M$. It is clear that $\rho_M :
M \rightarrow M \tensor{A} \coring{C}$ is a morphism of left
$\coring{C}'$--comodules if and only if $\lambda_M : M \rightarrow
\coring{C}' \tensor{A'} M$ is a morphism of right
$\coring{C}$--comodules. In this case, we say that $M$ is a
$\coring{C}'-\coring{C}$--bicomodule.

For any right $A$--module $X$, we will denote its right dual by
$X^* = \hom{A}{X}{A}$, which is a left $A$--module in a canonical
way. We will use the analogous notation ${}^*Y$ for left
$A$--modules $Y$. There is a canonical isomorphism
$\rend{\coring{C}}{\coring{C}} \cong \coring{C}^*$ that maps an
endomorphism $f$ onto $\epsilon \circ f$. The structure of ring of
$\rend{\coring{C}}{\coring{C}}$ is automatically transferred to
the convolution product $\coring{C}^*$ as defined in
\cite[Proposition 3.2]{Sweedler:1975}. Analogously, there is an
isomorphimsm of rings $\lend{\coring{C}}{\coring{C}}^{op} \cong
{}^*\coring{C}$. The coring $\coring{C}$ becomes a
${}^*\coring{C}-\coring{C}^*$--bimodule.

\section{Comatrix corings and the canonical map.}\label{comatrixcan}

Let $A, B$ be rings. A comatrix $A$--coring will be built on every
$B-A$--bimodule which is finitely generated and projective as a
right $A$--module. Sweedler's canonical corings and dual corings
are examples of comatrix corings. When the bimodule enjoys a
structure of right module over a given coring $\coring{C}$, a
canonical homomorphism of $A$--corings is shown to relate the
comatrix $A$--coring and the coring $\coring{C}$.

 Let ${}_B\Sigma_A$ be a
$B-A$--bimodule. Then $\Sigma^* = \hom{A}{\Sigma}{A}$ is
canonically endowed with a structure of $A-B$--bimodule, and
$\Sigma^* \tensor{B} \Sigma$ is an $A$--bimodule in a natural way.
Assume $\Sigma_A$ to be finitely generated and projective, and let
$\{ e_i^*,e_i \} \subseteq \Sigma^* \times \Sigma$ be a finite
dual basis, that is, the equality
\begin{equation}\label{dualbasis}
u = \sum_i e_i e_i^*(u) \qquad \forall u \in \Sigma
\end{equation}
holds. We can built on $\rcomatrix{B}{\Sigma}$ a canonical
structure of $A$--coring.

\begin{proposition}\label{comatrix}
The $A$--bimodule $\rcomatrix{B}{\Sigma}$ is an $A$--coring with
comultiplication
\[ \Delta : \rcomatrix{B}{\Sigma} \rightarrow
\rcomatrix{B}{\Sigma} \tensor{A} \rcomatrix{B}{\Sigma}
\]
defined by $\Delta(\varphi \tensor{B} u) = \sum_{i} \varphi
\tensor{B} e_i \tensor{A} e_i^* \tensor{B} u$ and counit
\[
\epsilon : \rcomatrix{B}{\Sigma} \rightarrow A
\]
 given by
$\epsilon(\varphi \tensor{B} u) = \varphi(u)$.  Moreover, there is
a ring anti-isomorphism ${}^*(\Sigma^* \tensor{B} \Sigma) \cong
\mathrm{End}({}_B \Sigma)$, where the first of these rings is the
left dual of $\Sigma^* \tensor{B} \Sigma$, endowed with the
convolution product.
\end{proposition}
\begin{proof}
First, we should check that $\Delta$ is well-defined. This
requires to prove that for every $b \in B$ and every pair
$(\varphi , u) \in \Sigma^* \times \Sigma$ one has
\[
\sum_i \varphi b \tensor{B} e_i \tensor{A} e_i^* \tensor{B} u =
\sum_i \varphi \tensor{B} e_i \tensor{A} e_i^* \tensor{B} bu
\]
Clearly, it suffices to show that $\sum_i b e_i \tensor{A} e_i^* =
\sum_i e_i \tensor{A} e_i^*b$. Using \eqref{dualbasis}, we compute
\begin{multline}\label{Blin}
\sum_i be_i \tensor{A} e_i^* = \sum_{i,k} e_ke_k^*(be_i)
\tensor{A} e_i^* = \sum_{i,k} e_k \tensor{A} e_k^*(be_i)e_i^* = \\
\sum_k e_k \tensor{A} (\sum_i (e_k^*b)(e_i)e_i^*) = \sum_k e_k
\tensor{A} e_k^*b,
\end{multline}
as desired. It is now routine to check that $\Delta$ and
$\epsilon$ are homomorphisms of $A$--bimodules. An easy
computation gives the coassociative and counitary properties. The
reader should realize that \eqref{dualbasis} is used again in the
proof of the counitary property. Finally, let us prove the stated
ring isomorphism. The isomorphism is given, at the level of
$K$--modules, by the composition
\[
{}^*(\Sigma^* \tensor{B} \Sigma) = \hom{A}{\Sigma^* \tensor{B}
\Sigma}{{}_AA} \cong \hom{B}{\Sigma}{{}^*(\Sigma^*)} \cong
\hom{B}{\Sigma}{\Sigma},
\]
where we have used one adjointness isomorphism and the canonical
isomorphism $\Sigma \cong {}^*(\Sigma^*)$. By computing explicitly
this composition is given by the assignment $f \mapsto \hat{f}$,
where $\hat{f}: \Sigma \rightarrow \Sigma$ is given by
\begin{equation}\label{hat}
\hat{f}(u) = \sum_ie_if(e_i^* \tensor{B} u), \text{ for every } f
\in {}^*(\Sigma^* \tensor{B} \Sigma).
\end{equation}
We are now ready to check that $\widehat{(-)}$ is a ring
anti-homomorphism. First, we have
$\hat{\epsilon}(u)=\sum_ie_i\epsilon(e_i^* \tensor{B}
u)=\sum_ie_ie_i^*(u)=u$, for every $u \in \Sigma$. Given, $f, g
\in {}^*(\Sigma^* \tensor{B} \Sigma)$, the convolution product
reads
\[
(f * g)(\varphi \tensor{B} u) = f(I \tensor{A} g)\Delta(\varphi
\tensor{B} u) = \sum_i f(\varphi \tensor{B} e_i g(e_i^* \tensor{B}
u)).
\]
Therefore,
\begin{multline*}
\widehat{f * g}(u) = \sum_j e_j(f*g)(e_j^* \tensor{B} u) =
\sum_{i,j}e_jf(e_j^* \tensor{B} e_i g(e_i^* \tensor{B} u)) = \\
\hat{f}(\sum_i e_i g(e_i^* \tensor{B} u)) = \hat{f}(\hat{g}(u)),
\end{multline*}
and we get that $\widehat{f*g} = \hat{f} \circ \hat{g}$. Since the
product in $\mathrm{End}({}_B\Sigma)$ is the opposite of the
composition, we have already proved that our map is an
anti-isomorphism of rings.
\end{proof}

\begin{remark}
Let $\{e_i^*,e_i\}_{1 \leq i \leq n} $ and $\{f_j^*,f_j\}_{1 \leq
j \leq m}$ be two dual bases for $\Sigma_A$. The well-known
isomorphism of $\rend{A}{\Sigma}$--bimodules
\begin{equation}\label{xi}
\xymatrix@R=0pt@C=30pt{\Sigma \tensor{A} \Sigma^*
\ar@{->}^{\xi}[r] & \rend{A}{\Sigma} \\ u \tensor{A} \varphi
\ar@{|->}[r] & [v \mapsto u\varphi(v) ] }
\end{equation}
easily shows that $\sum_{1 \leq i \leq n}e_i \tensor{A} e_i^* =
\sum_{1 \leq j \leq m}f_j \tensor{A} f_j^*$ and, hence, the
comultiplication of the comatrix $A$--coring
$\rcomatrix{B}{\Sigma}$ does not depend on the choice of the dual
basis of $\Sigma_A$.
\end{remark}

Our comatrix corings generalize two fundamental classes of
corings.

\begin{example}[Sweedler's canonical coring]
Let $B \rightarrow A$ a ring homomorphism. The bimodule ${}_BA_A$
is projective and finitely generated as a right $A$--module. The
corresponding comatrix $A$--coring is then isomorphic to the
\emph{canonical Sweedler's coring $A \tensor{B} A$} \cite[Example
1.2]{Sweedler:1975}. The comultiplication sends $a \tensor{B} a'$
onto $a \tensor{B} 1 \tensor{A} 1 \tensor{B} a'$, and the counit
is given by the multiplication of $A$.
\end{example}

\begin{example}[Dual coring]\label{dualcoring}
Let $A \rightarrow B$ a ring homomorphism. Assume that $B_A$ is
finitely generated and projective. Then, taking the dual with
respect to $A$, we have $B^* \tensor{B} B \cong B^*$, and this
isomorphism of $A$--bimodules becomes an isomorphism of
$A$--corings whenever we consider the $A$--coring structure on
$B^*$ obtained from the $A$--ring structure of $B$
\cite[3.7]{Sweedler:1975}.
\end{example}

A relevant feature of the $A$--coring $\Sigma^* \tensor{B} \Sigma$
is that the right $A$--module $\Sigma$ becomes a right $\Sigma^*
\tensor{B} \Sigma$--comodule in a canonical way. Its coaction is
defined as
\[
\xymatrix{\Sigma \ar^-{\rho_{\Sigma}}[r] & \Sigma \tensor{A}
\Sigma^* \tensor{B} \Sigma, & & (u \ar@{|->}[r] & \sum_i e_i
\tensor{A} e_i^* \tensor{B} u)},
\]
which is clearly a left $B$--linear map, in other words, $\Sigma$
is $B-(\Sigma^* \tensor{B} \Sigma)$--bicomodule. This comodule
plays a relevant role.

\begin{proposition}\label{subgenra}
The coring $\Sigma^* \tensor{B} \Sigma$ is, as a right comodule,
generated by $\Sigma$. Therefore, every right
$\rcomatrix{B}{\Sigma}$--comodule is isomorphic to a subcomodule
of a quotient of $\Sigma^{(I)}$, for a suitable index set $I$.
Moreover,
\[
\rend{\rcomatrix{B}{\Sigma}}{\Sigma} = \{ f \in
\rend{A}{\Sigma}|\,\, f \tensor{B} x = 1 \tensor{B} f(x), \text{
for every } x \in \Sigma \}
\]
 and, in particular, the canonical
ring homomorphism $B \rightarrow \rend{A}{\Sigma}$ factorizes
throughout $\rend{\rcomatrix{B}{\Sigma}}{\Sigma}$.
\end{proposition}
\begin{proof}
For the first statement, it is enough to prove that every
generator $\varphi \tensor{B} u \in \Sigma^* \tensor{B} \Sigma$
belongs to the image of a morphism of right $\Sigma^* \tensor{B}
\Sigma$--comodules $f : \Sigma \rightarrow \rcomatrix{B}{\Sigma}$.
This is fulfilled by the map defined by $f(u) = \varphi \tensor{B}
u$, which is easily proved to be a homomorphism of comodules. For
the second statement, let $\rho_M : M \rightarrow M \tensor{A}
\rcomatrix{B}{\Sigma}$ be a right comodule. The own structure map
$\rho_M$ is a morphism of comodules which splits as a right
$A$--module map. Therefore, $M$ is isomorphic to a subcomodule of
$M \tensor{A} \rcomatrix{B}{\Sigma}$, which is now easily shown to
be a quotient of a coproduct of copies of $\Sigma$. The second
statement follows from a straightforward computation.
\end{proof}

\begin{remark}\label{comizq}
Of course $\Sigma^*$ is a left $\rcomatrix{B}{\Sigma}$--comodule
with a right $B$--linear coaction
\[
\xymatrix{ \Sigma^* \ar^-{\lambda_{\Sigma^*}}[r] &
\rcomatrix{B}{\Sigma} \tensor{A} \Sigma^*, & & (\varphi \mapsto
\sum_i \varphi \tensor{B} e_i \tensor{A} e_i^*)}.
\]
Moreover, ${}_{\rcomatrix{B}{\Sigma}}\Sigma^*$ satisfies the left
version of Proposition \ref{subgenra}, and the right convolution
ring $(\rcomatrix{B}{\Sigma})^*$ is a ring anti-isomorphic to
$\rend{B}{\Sigma^*}$.
\end{remark}

Now, let $\coring{C}$ be any $A$--coring, and assume $\Sigma$ to
be a right $\coring{C}$--comodule with coaction $\rho_{\Sigma} :
\Sigma \rightarrow \Sigma \tensor{A} \coring{C}$.  From now on, we
will denote $S = \rend{A}{\Sigma}$ and $T =
\rend{\coring{C}}{\Sigma}$. Then $\Sigma$ becomes an
$S-A$--bimodule and $\rho_{\Sigma}$ is a homomorphism of
$S-A$--bimodules. We keep in mind that $T$ is a subring of $S$.

\begin{proposition}\label{pi}
If $\Sigma_{\coring{C}}$ is a comodule such that $\Sigma_A$ is
finitely generated and projective, then the map $can :
\rcomatrix{T}{\Sigma} \rightarrow \coring{C}$ defined as the
composition
\[
\xymatrix {\rcomatrix{T}{\Sigma} \ar^-{\Sigma^* \tensor{B}
\rho_{\Sigma}}[rr] & & \rcomatrix{T}{\Sigma} \tensor{A} \coring{C}
\ar^-{\epsilon \tensor{A} \coring{C}}[rr] & & A \tensor{A}
\coring{C} \cong \coring{C} }
\]
is a homomorphism of $A$--corings.
\end{proposition}
\begin{proof}
By construction $can$ is $A$--bilinear. We need to check the
identities
\begin{equation}\label{coringmorf}
\Delta_{\coring{C}} \circ can = (can \tensor{A} can) \circ \Delta
\qquad \textrm{and} \qquad \epsilon_{\coring{C}} \circ can =
\epsilon.
\end{equation}
If $\{ e_i^*, e_i \} \subseteq \Sigma^* \times \Sigma$ is a dual
basis, then $\rho_{\Sigma}(e_j) = \sum_i e_i \tensor{A} \rho_{ij}$
for each $j$. Then the first identity in \eqref{coringmorf} is
equivalent to
\begin{equation}\label{coringmorf2}
\sum_k e_i^*(e_k)\Delta_{\coring{C}}(\rho_{kj}) =
\sum_{k,l,m}e_i^*(e_l)\rho_{lk} \tensor{A} e_k^*(e_m)\rho_{mj}
\qquad \forall i, j
\end{equation}
The fact that $\Sigma$ is a right $\coring{C}$--comodule gives
\[
\sum_k e_k \tensor{A} \Delta_{\coring{C}}(\rho_{kj}) = \sum_{k,l}
e_k \tensor{A} \rho_{kl} \tensor{A} \rho_{lj} \qquad \forall j
\]
Therefore,
\begin{equation}\label{unaec}
\sum_k e_i^*(e_k)\Delta_{\coring{C}}(\rho_{kj}) = \sum_{k,l}
e_i^*(e_k) \rho_{kl} \tensor{A} \rho_{lj} \qquad \forall i, j
\end{equation}
On the other hand,
\[
\rho_{\Sigma}(e_m) = \rho_{\Sigma}(\sum_k e_ke_k^*(e_m)) = \sum_k
\rho_{\Sigma}(e_k)e_k^*(e_m) = \sum_{l,k} e_l \tensor{A}
\rho_{lk}e_k^*(e_m) \qquad \forall m;
\]
thus
\[
\sum_l e_l \tensor{A} \rho_{lm} = \sum_{l,k}e_l \tensor{A}
\rho_{lk}e_k^*(e_m) \quad \textrm{and} \quad \sum_l e_i^*(e_l)
\rho_{lm} = \sum_{l,k}e_i^*(e_l) \rho_{lk}e_k^*(e_m) \quad \forall
i,m
\]
Finally, we get
\[
\sum_{k,l,m}e_i^*(e_l)\rho_{lk} \tensor{A} e_k^*(e_m)\rho_{mj} =
\sum_{k,l,m}e_i^*(e_l)\rho_{lk}e_k^*(e_m) \tensor{A} \rho_{mj} =
\sum_{l,m}e_i^*(e_l)\rho_{lm} \tensor{A} \rho_{mj} \quad \forall
i, j
\]
This implies, in view of \eqref{unaec}, the identity
\eqref{coringmorf2}.\\ The second identity in \eqref{coringmorf}
holds since, for every $i, j$, we have
\begin{multline*}
\epsilon_{\coring{C}}(can(e_i^* \tensor{} e_j)) =
\epsilon_{\coring{C}}(\sum_k e_i^*(e_k)\rho_{kj}) = \sum_k
e_i^*(e_k)\epsilon_{\coring{C}}(\rho_{kj}) \\ = e_i^*(\sum_k e_k
\epsilon_{\coring{C}}(\rho_{kj}) ) = e_i^*(e_j) = \epsilon(e_i^*
\tensor{} e_j)
\end{multline*}
\end{proof}

The following examples suggest that the $can$ map defined in
Proposition \ref{pi} is an interesting object for research.
Moreover, it generalizes canonical maps previously considered in
the theories of Hopf modules and noncommutative Galois extensions.

\begin{example}\label{canescan}
Let us assume that our $A$--coring $\coring{C}$ has a grouplike
element $g$, which is equivalent, by \cite[Lemma
5.1]{Brzezinski:2000unp}, to endow $A$ with a structure of right
comodule over $\coring{C}$. In this case $T =
\rend{\coring{C}}{A}$ is nothing but the subring of coinvariants
\cite[Proposition 2.2]{Brzezinski:2001unp} of $A$. In this case,
the map $can : A \tensor{T} A \rightarrow \coring{C}$ is
determined by the condition $can(1 \tensor{T} 1) = g$. Therefore,
the coring $\coring{C}$ is Galois in the sense of \cite[Definition
5.3]{Brzezinski:2000unp} if and only if $can$ is an isomorphism.
It is convenient to point out here that our homomorphism $can$
generalizes the original canonical map considered by
\cite{Schneider:1990a} and, in fact, the map $can$ defined in
\cite[Definition 2.1]{Brzezinski:1999}.
\end{example}

\begin{example}
Let $G$ be a finite group of ring automorphims of $A$, and let $R
= G*A$ be associated crossed product. The ring $A$ embeds
canonically in $R$ and, by construction, $R_A$ is free with basis
$G$ and we can consider the corresponding comatrix coring $R^*
\tensor{R} R \cong R^*$ (see Example \ref{dualcoring}). Let us
show that the ``trace map'' $g : R \rightarrow A$ defined by
$g(\sum_{\sigma \in G}\sigma a_{\sigma}) = \sum_{\sigma \in
G}a_{\sigma}$ is a grouplike element for $R^*$.  Accordingly with
\cite[Theorem 3]{Kleiner:1984}, we need just to check that $g$
acts as the identity on $A$, which is obviously the case, and that
$Ker g$ is a right ideal of $R$. This last condition can be
checked in a straightforward way taking that the trace map is
invariant under translations into account. Thus, $A$ is a right
$R^*$--comodule and we have the homomorphism of $A$--corings $can
: A \tensor{T} A \rightarrow R^*$ determined by $can(1 \tensor{T}
1) = g$, where $T$ is the subring of $g$--coinvariants of $A$. An
easy computation shows that $T$ is already the subring of
$G$--invariants of $A$. Now, the composite homomorphism of rings
\[
\xymatrix{R \cong {}^*(R^*) \ar^-{{}^*can}[rr] & & {}^*(A
\tensor{T} A) \cong \lend{T}{A}}
\]
is precisely the map $\delta$ defined in \cite{Kanzaki:1964} (or
$j$ in \cite{DeMeyer:1965}). There,  the extension $T \subseteq A$
is said to be $G$--Galois whenever $\delta$ is an isomorphism and
${}_TA$ is finitely generated and projective. Of course, $\delta$
is an isomorphism if and only if $can$ is an isomorphism.
\end{example}

The homomomorphism of $A$--corings $can : \rcomatrix{T}{\Sigma}
\rightarrow \coring{C}$ leads to the functor
\[
CAN : \rcomod{\rcomatrix{T}{\Sigma}} \rightarrow
\rcomod{\coring{C}}
\]
 which sends a comodule $\rho_M : M \rightarrow
M \tensor{A} \rcomatrix{T}{\Sigma}$ onto the
$\coring{C}$--comodule
\[
\xymatrix{M \ar^-{\rho_M}[r] & M \tensor{A} \rcomatrix{T}{\Sigma}
\ar^-{M \tensfun{A} can}[rr] & & M \tensor{A} \coring{C}}
\]
This is an example of induction functor (see
\cite[5.2]{Gomez:2002} for details).

\begin{proposition}\label{CAN}
If $\Sigma_{\coring{C}}$ is a comodule such that $\Sigma_A$ is
finitely generated and projective,  then $T =
\rend{\rcomatrix{T}{\Sigma}}{\Sigma}$ and we have a commuting
diagram of functors
\[
\xymatrix{\rcomod{\rcomatrix{T}{\Sigma}} \ar^-{CAN}[rr]  & & \rcomod{\coring{C}} \\
 & \rmod{T} \ar^{-
\tensfun{T} \Sigma }[ul] \ar_{- \tensfun{T} \Sigma}[ur]&}
\]
\end{proposition}
\begin{proof}
 By
Proposition \ref{subgenra}, $T \subseteq
\rend{\rcomatrix{T}{\Sigma}}{\Sigma}$. Conversely, let $f \in
\mathrm{End}(\Sigma_{\rcomatrix{T}{\Sigma}})$; the following
diagram is clearly commutative
\[
\xymatrix@R=30pt@C=50pt{ \Sigma \ar@{->}^-{\rho_{\Sigma}}[r]
\ar@{->}_-{f}[d] & \Sigma \tensor{A} \rcomatrix{T}{\Sigma}
\ar@{->}^-{f \tensfun{A} \rcomatrix{T}{\Sigma}}[d]
\ar@{->}^-{\Sigma \tensfun{A} can}[r] & \Sigma
\tensor{A} \coring{C} \ar@{->}^-{f \tensfun{A} \coring{C}}[d] \\
\Sigma \ar@{->}^-{\rho_{\Sigma}}[r] & \Sigma \tensor{A}
\rcomatrix{T}{\Sigma} \ar@{->}^-{\Sigma \tensfun{A} can}[r] &
\Sigma \tensor{A} \coring{C}. }
\]
Now, an easy computation shows that $(\Sigma \tensor{A}
can)\rho_{\Sigma}$ is just the structure map for
$\Sigma_{\coring{C}}$. Thus $f$ is right $\coring{C}$--colinear,
that is, $f \in T$. Observe that we have already shown that
$CAN(\Sigma_{\rcomatrix{T}{\Sigma}}) = \Sigma_{\coring{C}}$. This
implies that $CAN((X \tensor{T} \Sigma)_{\rcomatrix{T}{\Sigma}}) =
(X \tensor{T}\Sigma)_{\coring{C}}$ for every $X \in \rmod{T}$, as
desired.
\end{proof}

\section{Corings with a finitely generated and projective generator}\label{gl-dethy}

We give a complete description in terms of comatrix corings of
corings having a finitely generated projective generator.
Furthermore, our result generalizes \cite[Theorem
5.6]{Brzezinski:2000unp} to corings which possibly have not
grouplike elements and, therefore, it is ultimately a
generalization of \cite[Theorem 1]{Schneider:1990a} and
\cite[Teorema]{Cipolla:1976}.

Let $\Sigma$ be a right comodule over an $A$--coring $\coring{C}$,
and let $T = \rend{\coring{C}}{\Sigma}$ its endomorphism ring. The
structure map of $\Sigma$ is $T$--linear and, thus, we have the
functor $- \tensor{T} \Sigma : \rmod{T} \rightarrow
\rcomod{\coring{C}}$. Recall that $\hom{\coring{C}}{\Sigma}{-} :
\rcomod{\coring{C}} \rightarrow \rmod{T}$ is right adjoint to $-
\tensor{T} \Sigma$. Let $\chi : \hom{\coring{C}}{\Sigma}{-}
\tensor{T} \Sigma \rightarrow 1$ the counit of this adjunction.

The proof of our main theorem will be better understood if we
isolate a technical fact which, in view of \cite[Proposition
1.1]{Caenepeel/Vercruysse/Wang:2002unp} and \cite[Theorem
2.2]{Abuhlail:2002unp} seems to be of independent interest.

\begin{lemma}\label{canchi}
Let $\Sigma$ be a right $\coring{C}$--comodule such that
$\Sigma_A$ is finitely generated and projective. The counit
$\chi_{\coring{C}}$ at $\coring{C}$
 is an isomorphism if and only if $can$ is an
isomorphism.
\end{lemma}
\begin{proof}
Making use of the isomorphism
$\hom{\coring{C}}{\Sigma}{\coring{C}} \cong \Sigma^*$, we have
that $can : \rcomatrix{T}{\Sigma} \rightarrow \coring{C}$ can be
written as the composite
\[
\xymatrix{\rcomatrix{T}{\Sigma} \cong
\hom{\coring{C}}{\Sigma}{\coring{C}} \tensor{T} \Sigma
\ar^-{\chi_{\coring{C}}}[r] & \coring{C}. }
\]
\end{proof}

\begin{theorem}\label{generador}
Let $\coring{C}$ be an $A$--coring, and $\Sigma_{\coring{C}}$ a
right $\coring{C}$--comodule. Consider the ring extension $T
\subseteq S$, where $T = \mathrm{End}(\Sigma_{\coring{C}})$ and $S
= \rend{A}{\Sigma}$. The following statements are equivalent
\begin{enumerate}[(i)]
\item ${}_A\coring{C}$ is flat and $\Sigma_{\coring{C}}$ is a finitely generated and projective
generator for $\rcomod{\coring{C}}$;
\item  ${}_A\coring{C}$ is flat, $\Sigma_A$ is a
finitely generated and projective $A$--module, $can:
\rcomatrix{T}{\Sigma} \rightarrow \coring{C}$ is an isomorphism of
$A$--corings, and $- \tensor{T} \Sigma : \rmod{T} \rightarrow
\rcomod{\rcomatrix{T}{\Sigma}}$ is an equivalence of categories;
\item $\Sigma_A$ is a
finitely generated and projective $A$--module, $can:
\rcomatrix{T}{\Sigma} \rightarrow \coring{C}$ is an isomorphism of
$A$--corings, and ${}_{T} \Sigma$ is a faithfully flat module.
\item  ${}_A\coring{C}$ is flat, $\Sigma_A$ is a
finitely generated and projective $A$--module, $can:
\rcomatrix{T}{\Sigma} \rightarrow \coring{C}$ is an isomorphism of
$A$--corings, and ${}_TS$ is faithfully flat.
\end{enumerate}
\end{theorem}
\begin{proof}
$(i) \Rightarrow (ii)$ Since ${}_A\coring{C}$ is flat, it follows
from \cite[Proposition 1.2]{ElKaoutit/Gomez/Lobillo:2001pr} that
$\rcomod{\coring{C}}$ is a Grothendieck category and the forgetful
functor $U : \rcomod{\coring{C}} \rightarrow \rmod{A}$ is exact.
Moreover, it has an exact right adjoint $- \tensor{A} \coring{C}$.
This implies that $\Sigma_A$ is finitely generated and projective.
Recall that $\hom{\coring{C}}{\Sigma}{-} : \rcomod{\coring{C}}
\rightarrow \rmod{T}$ is right adjoint to $ - \tensor{T} \Sigma$
and, since $\Sigma_{\coring{C}}$ is a finitely generated and
projective generator, it is already an equivalence of categories.
In particular, the counity of the adjunction $\chi :
\hom{\coring{C}}{\Sigma}{-} \tensor{T} \Sigma \rightarrow 1$ is a
natural isomorphism.  By Lemma \ref{canchi}, $CAN :
\rcomod{\rcomatrix{T}{\Sigma}} \rightarrow \rcomod{\coring{C}}$ is
an equivalence of categories. By Proposition \ref{CAN}, we have
that $- \tensor{T} \Sigma : \rmod{T} \rightarrow
\rcomod{\rcomatrix{T}{\Sigma}}$ is an equivalence of categories.
\\
$(ii) \Rightarrow (iii)$ The functor $- \tensor{T} \Sigma :
\rmod{T} \rightarrow \rcomod{\rcomatrix{T}{\Sigma}}$ is obviously
faithful and exact. Since $\rcomatrix{T}{\Sigma} \cong \coring{C}$
is flat as a left $A$--module, we have, by \cite[Proposition 1.2
]{ElKaoutit/Gomez/Lobillo:2001pr}, that the forgetful functor $U :
\rcomod{\rcomatrix{T}{\Sigma}} \rightarrow \rmod{A}$ is faithful
and exact. Therefore, the functor $- \tensor{T} \Sigma : \rmod{T}
\rightarrow \rmod{A}$ is faithful and exact, that is,
${}_T\Sigma$ is a faithfully flat module.\\
$(iii) \Rightarrow (i)$. Let $\Omega = \rcomatrix{T}{\Sigma}$,
which is flat as a left $A$--module because ${}_T\Sigma$ is flat.
The isomorphism of $A$--corings $can:\Omega\cong \coring{C}$ gives
then that ${}_A\coring{C}$ is flat. Consider on $\Sigma^*$ the
left $\Omega$--comodule structure given in Remark \ref{comizq}. A
straightforward computation shows that $\Sigma^*$ is an
$\Omega-T$--bicomodule, where $T$ acts canonically on $\Sigma^*$.
The functor $ - \tensor{A} \Sigma^* : \rmod{A} \rightarrow
\rmod{T}$ is right adjoint to $- \tensor{T} \Sigma: \rmod{T}
\rightarrow \rmod{A}$. By \cite[Proposition 4.2]{Gomez:2002}, the
cotensor product functor $- \cotensor{\Omega}\Sigma^* :
\rcomod{\Omega} \rightarrow \rmod{T}$ is right adjoint to $-
\tensor{T} \Sigma : \rmod{T} \rightarrow \rcomod{\Omega}$.  Since
${}_T\Sigma$ is flat we have, by \cite[Lemma 2.2]{Gomez:2002}, the
isomorphism $(M \cotensor{\Omega} \Sigma^*)\tensor{T} \Sigma \cong
M \cotensor{\Omega}(\Sigma^* \tensor{T} \Sigma) = M
\cotensor{\Omega} \Omega \cong M$, which turns out to be inverse
to the counity of the adjunction at $M \in \rcomod{\Omega}$.
Moreover, if $\eta_X : X \rightarrow (X \tensor{T}
\Sigma)\cotensor{\Omega}\Sigma^*$ is the unity of the adjunction
at $X \in \rmod{T}$, then an inverse to $\eta_X \tensor{T} \Sigma$
is obtained by the isomorphism $((X \tensor{T} \Sigma)
\cotensor{\Omega} \Sigma^*))\tensor{T} \Sigma \cong (X \tensor{T}
\Sigma) \cotensor{\Omega} (\Sigma^* \tensor{T} \Sigma) \cong X
\tensor{T} \Sigma$. Since ${}_T\Sigma$ is faithful, we get that
$\eta_X$ is an isomorphism and, hence, $- \tensor{T} \Sigma
:\rmod{T} \rightarrow \rcomod{\Omega}$ is an equivalence of
categories with inverse $- \cotensor{\Omega} \Sigma^*$. It follows
from Proposition \ref{CAN} that $- \tensor{T} \Sigma : \rmod{T}
\rightarrow \rcomod{\coring{C}}$ is an equivalence, as $can :
\Omega \rightarrow \coring{C}$ is an isomorphism. Therefore,
$\Sigma_{\coring{C}}$ is a finitely generated and projective
generator for $\rcomod{\coring{C}}$.\\
$(iii) \Rightarrow (iv)$. The exact and faithful functor $-
\tensor{T} \Sigma : \rmod{T} \rightarrow \rmod{A}$ decomposes as
$- \tensor{T} \Sigma \simeq (- \tensor{S} \Sigma)\circ(-
\tensor{T} S)$. This implies that $- \tensor{T} S$ is an exact and
faithful functor, as $- \tensor{S} \Sigma$ is faithful (because
$\Sigma_A$ is finitely generated and projective). Therefore, ${}_TS$ is faithfully flat. \\
$(iv) \Rightarrow (iii)$. It suffices to show that ${}_T\Sigma$ is
faithfully flat. Consider a short exact sequence in $\rmod{T}$,
$\xymatrix{0 \ar@{->}[r]& Y \ar@{->}^-{f}[r] & Y' \ar@{->}^-{g}[r]
& Y'' \ar@{->}[r] & 0}$, and let $Z$ denote the kernel of the
morphism $f\tensor{T}\Sigma$ in the category
$\rcomod{\rcomatrix{T}{\Sigma}}$. Since the forgetful functor
$\rcomod{\rcomatrix{T}{\Sigma}} \rightarrow \rmod{A}$ is exact,
this kernel coincides with the kernel computed in $\rmod{A}$. We
thus get a commutative diagram with exact rows
\[
\xymatrix{0 \ar@{->}[r] & Z \tensor{A} \Sigma^* \ar@{->}[r] & Y
\tensor{T} \Sigma \tensor{A} \Sigma^* \ar@{->}[r]
\ar@{->}^{\simeq}[d] & Y' \tensor{T} \Sigma \tensor{A} \Sigma^*
\ar@{->}[r] \ar@{->}^-{\simeq}[d] & Y'' \tensor{T} \Sigma
\tensor{A} \Sigma^* \ar@{->}[r] \ar@{->}^-{\simeq}[d] & 0 \\ & 0
\ar@{->}[r] & Y \tensor{T} S \ar@{->}[r] & Y'\tensor{T} S
\ar@{->}[r] & Y'' \tensor{T} S \ar@{->}[r] & 0. }
\]
Therefore, $Z \tensor{A} \Sigma^*=0$, which implies that $Z
\tensor{A} \Sigma^* \tensor{T} \Sigma =0$, and thus $Z=0$, since
$Z \in \rcomod{\rcomatrix{T}{\Sigma}}$. Therefore, ${}_T\Sigma$ is
flat. Now, given any right $T$--module $Y$ such that $Y \tensor{T}
\Sigma = 0$, we have $Y \tensor{T} \Sigma \tensor{A} \Sigma^*
\cong Y \tensor{T} S = 0$, and so $Y=0$.  Thus, ${}_T\Sigma$ is a
faithfully flat module, and this finishes the proof.
\end{proof}

\begin{remark}
The flatness of ${}_A\coring{C}$ cannot be dropped in the
statements $(ii)$ and $(iv)$. Counterexamples can be obtained as
follows. Let $e \in A$ an idempotent and write $f = 1-e$. Assume
that $fAe = 0$ and let $I = eA$. Then $I$ becomes an $A$--coring
whose comultiplication is given by the canonical isomorphism $I
\cong I \tensor{A} I$ and the counit is just the inclusion $I
\subseteq A$. In this case, $\rend{I}{I} = \rend{A}{I} \cong eAe$,
and, hence, $T = S$. Moreover, this easily implies that $can : I^*
\tensor{eAe} I \cong Ae \tensor{eAe} I \cong I$. The right
$I$--comodules are those right $A$--modules satisfying that the
canonical map $M \tensor{A} I \rightarrow M$ is an isomorphism.
This allows to prove that $- \tensor{eAe} I :\rmod{eAe}
\rightarrow \rcomod{I}$ is an equivalence of categories. However,
${}_AI$ is not flat unless ${}_{eAe}I$ is.
\end{remark}

Next, we shall locate Theorem \ref{generador} within the recent
developments Galois corings with grouplike elements.

\begin{definition}
Let $\coring{C}$ be an $A$--coring with a right
$\coring{C}$--comodule $\Sigma$ such that $\Sigma_A$ is finitely
generated and projective. The coring will be said to be
\emph{Galois} if $can : \rcomatrix{T}{\Sigma} \rightarrow
\coring{C}$ is an isomorphism, where $T =
\rend{\coring{C}}{\Sigma}$. When $\Sigma = A$, this definition
coincides with the given by T. Brzezi\'{n}ski for corings with a
grouplike element \cite[Definition 5.3]{Brzezinski:2000unp}. In
view of Theorem \ref{generador}, a ring extension of the form $T
\rightarrow \rend{A}{\Sigma}$ could be called a
\emph{$\coring{C}$--Galois ring extension} whenever $\Sigma$ is a
right $\coring{C}$--comodule such that $\Sigma_A$ is finitely
generated and projective and $T = \rend{\coring{C}}{\Sigma}$.
Coring Galois extensions in \cite{Brzezinski:2000unp} or
\cite{Caenepeel/Vercruysse/Wang:2002unp} are then obtained with
$\Sigma = A$.
\end{definition}

We easily get now.

\begin{corollary}\cite[Theorem 5.6]{Brzezinski:2000unp}
Let $\coring{C}$ be an $A$--coring with a grouplike element, and
$T$ be the subring of coinvariant elements of $A$. If $\coring{C}$
is Galois and ${}_TA$ is faithfully flat, then $- \tensor{T} A :
\rmod{T} \rightarrow \rcomod{\coring{C}}$ is an equivalence of
categories. Conversely, if $- \tensor{T} A$ is an equivalence of
categories, then $\coring{C}$ is Galois. In this case if
${}_A\coring{C}$ is flat, then ${}_TA$ is faithfully flat.
\end{corollary}
\begin{proof}
Put $\Sigma = A$. The corollary follows from Example
\ref{canescan}, Proposition \ref{CAN}, Lemma \ref{canchi} and,
mainly, Theorem \ref{generador}.
\end{proof}

Galois corings with grouplike element has been also recently
considered from the point of view of category equivalences by J.
Y. Abuhlail \cite{Abuhlail:2002unp}. Some of his results can be
easily derived from our set up.

\begin{remark}
In view of Lemma \ref{canchi}, $can$ is an isomorphism if and only
if $\chi_{\coring{C}}$ is an isomorphism. Taking $\Sigma = A$ in
Theorem \ref{generador} we obtain that for if $\coring{C}$ has a
grouplike element and ${}_A\coring{C}$ is flat then $A$ is a
projective generator for $\rcomod{\coring{C}}$ if and only if
${}_TA$ is faithfully flat and $\chi_{\coring{C}}$ is bijective.
This has been recently proved by J. Y. Abuhlail under the
additional condition ``${}_A\coring{C}$ is locally projective''
(see \cite[Theorem 2.4]{Abuhlail:2002unp}).
\end{remark}

\begin{remark}
It follows from Lemma \ref{canchi} and from the proof of $(iii)
\Rightarrow (i)$ in Theorem \ref{generador} that the counit
$\chi_M$ is an isomorphism for every $M \in \rcomod{\coring{C}}$
if and only if ${}_T\Sigma$ is flat and $\coring{C}$ is Galois.
Taking $\Sigma = A$, we get that the coring $\coring{C}$ with a
grouplike satisfies the ``Weak Structure Theorem'' if and only if
${}_TA$ is flat and $\coring{C}$ is Galois. This has been proved
in \cite[Theorem 2.2]{Abuhlail:2002unp} under the additional
condition ``${}_A\coring{C}$ is locally projective''.
\end{remark}

The interest in corings has been partly recovered because the
theory of entwined modules (and, henceforth, of Hopf modules) can
be subsumed in the theory of comodules over certain corings
\cite[Proposition 2.2]{Brzezinski:2000unp}.

\begin{example}
Let $(A,C)_{\psi}$ be an entwining structure over $K$ and assume
there is an entwined module $\Sigma$ such that $\Sigma_A$ is
finitely generated and projective. We have our canonical map $can
: \Sigma^* \tensor{T} \Sigma \rightarrow A \tensor{} C$, where $T$
is the ring of endomorphisms of $\Sigma$ as an entwined module. If
$can$ is bijective we have a special type of Galois $A$--corings
without grouplike elements. This process can be reversed: start
with a coalgebra $C$ and a finitely generated and projective right
module $\Sigma$ over an algebra $A$. Let $\rho : \Sigma
\rightarrow \Sigma \tensor{} C$ a structure of right $C$--comodule
over $\Sigma$ and define $T = \{ t \in \rend{A}{\Sigma} | \rho(tu)
= t\rho(u) \mbox{ for every } u \in \Sigma \}$ (that is, $T$ is
the ring of all endomorphisms of $\Sigma$ which are $A$--linear
and $C$--colinear). Then define $can(\varphi \tensor{T} u) = \sum
\varphi(u_{(0)}) \tensor{} u_{(1)}$. If this $can$ is bijective
(which could be the new more general definition of $C$--Galois
extension $T \subseteq \rend{A}{\Sigma}$) then there is a unique
entwining structure $(A,C)_{\psi}$ making $\Sigma$ a right
entwined module (to check this, first use propositions
\ref{comatrix} and  \ref{subgenra}  to transfer the structure of
$A$--coring of $\rcomatrix{T}{\Sigma}$ and of right
$\rcomatrix{T}{\Sigma}$--comodule of $\Sigma$, respectively; and
then \cite[Proposition 2.2]{Brzezinski:2000unp} to interpret
everything in terms of entwining structure). Taking $\Sigma = A$,
we obtain \cite[Theorem 2.7]{Brzezinski/Hajac:1999}.
\end{example}

The relationship between Noncommutative Descent Theory and Galois
corings with grouplike is known (see \cite{Brzezinski:2000unp} and
\cite{Caenepeel/Vercruysse/Wang:2002unp}). We will derive from our
analysis of Galois corings without grouplike elements a Descent
Theory for ring extensions of the form $B \rightarrow
\rend{A}{\Sigma}$, where $\Sigma_A$ is finitely generated and
projective. Of course, our sufficient conditions to have the
Descent Theorem are given on the bimodule $\Sigma$. Once again,
the case $\Sigma = A$ collapses with the classical theory.

\begin{lemma}\label{caniso}
Let ${}_B\Sigma_A$ be a $B-A$--bimodule such that $\Sigma_A$ is
finitely generated and projective, and let $T =
\rend{\rcomatrix{B}{\Sigma}}{\Sigma}$. Then the canonical map $can
: \rcomatrix{T}{\Sigma} \rightarrow \rcomatrix{B}{\Sigma}$ is an
isomorphism of $A$--corings.
\end{lemma}
\begin{proof}
Let $B \rightarrow T$ the homomorphism of rings given in
Proposition \ref{subgenra}. Denote by $\omega:
\rcomatrix{B}{\Sigma} \rightarrow \rcomatrix{T}{\Sigma}$ the
obvious map which sends $\varphi \tensor{B} x \mapsto \varphi
\tensor{T} x$. Let us check that $\omega$ is the inverse of $can$.
Given $\varphi \in \Sigma^*$ and $x \in \Sigma$, we have
$can(\omega(\varphi \tensor{B} x))= \sum_i \varphi(e_i)e_i^*
\tensor{B} x = \varphi \tensor{B} x$, and $\omega(can(\varphi
\tensor{T}x)= \omega(\sum_i\varphi(e_i)e_i^* \tensor{B} x)=
\omega(\varphi\tensor{B}x)=\varphi\tensor{T}x$. Hence, $can$ is an
isomorphism of $A$--corings.
\end{proof}

\begin{theorem}[Generalized Descent for Modules]\label{Descent}
Let ${}_B\Sigma_A$ be a $B-A$--bimodule such that $\Sigma_A$ is
finitely generated and projective. Then ${}_B\Sigma$ is faithfully
flat if and only if $- \tensor{B} \Sigma : \rmod{B} \rightarrow
\rcomod{\rcomatrix{B}{\Sigma}}$ is an equivalence of categories
and ${}_A(\rcomatrix{B}{\Sigma})$ is flat. In such a case, the
canonical map $\lambda : B \rightarrow
\rend{\rcomatrix{B}{\Sigma}}{\Sigma}$ is a ring isomorphism.
\end{theorem}
\begin{proof}
By Proposition \ref{subgenra}, we have that
\[
T = \rend{\rcomatrix{B}{\Sigma}}{\Sigma} =  \{ f \in
\rend{A}{\Sigma}|\,\, f \tensor{B} x = 1 \tensor{B} f(x), \text{
for every } x \in \Sigma \},
\]
which shows that the canonical map $B \tensor{B} \Sigma
\rightarrow T \tensor{B} \Sigma$ is an isomorphism. Therefore,
when ${}_B\Sigma$ is assumed to be flat, one deduces that $Ker
\lambda \tensor{B} \Sigma = coKer \lambda \tensor{B} \Sigma = 0$.
Thus, if ${}_B\Sigma$ is faithfully flat, then $\lambda$ is an
isomorphism of rings, and we can apply Lemma \ref{caniso} and
Theorem \ref{generador} to obtain that $ - \tensor{B} \Sigma:
\rmod{B} \rightarrow \rcomod{\rcomatrix{B}{\Sigma}}$ is an
equivalence of categories. Conversely, if we assume such an
equivalence, then $\Sigma$ is a finitely generated projective
generator for $\rcomod{\rcomatrix{B}{\Sigma}}$. We deduce from
Theorem \ref{generador} that $ - \tensor{T} : \rmod{T} \rightarrow
\rcomod{\rcomatrix{B}{\Sigma}}$ is an equivalence of categories
and ${}_T\Sigma$ is faithfully flat. Therefore, in the commutative
diagram of functors
\[
\xymatrix{\rmod{B} \ar^-{- \tensfun{B} \Sigma}[rr] & &
\rcomod{\rcomatrix{B}{\Sigma}} \\
\rmod{T} \ar^-{- \tensfun{T} \Sigma}[rr] \ar_-{F}[u] & &
\rcomod{\rcomatrix{T}{\Sigma}} \ar_-{CAN}[u]},
\]
where $F : \rmod{T} \rightarrow \rmod{B}$ is the restriction of
scalars functor associated to $\lambda : B \rightarrow T$, the
other three functors are equivalences of categories. This shows
that $\lambda$ is an isomorphism, which proves the Theorem.
\end{proof}

M. Cipolla \cite{Cipolla:1976} give a Descent Theorem for a
homomorphism of noncommutative rings $B \rightarrow A$. As T.
Brzezinski pointed out \cite{Brzezinski:2000unp}, the category of
descent data \cite{Nuss:1997} is precisely the category $\rcomod{A
\tensor{B} A}$ of right comodules over $A \tensor{B} A$. As a
consequence of Theorem \ref{Descent} we obtain Cipolla's main
result \cite[Teorema]{Cipolla:1976} (see also \cite[Theorem
3.8]{Nuss:1997}).

\begin{corollary}[Descent of Modules]\label{classicdescent}
Let $B \rightarrow A$ be a ring homomomorphism. If ${}_BA$ is
faithfully flat, then the functor $- \tensor{B} A : \rightarrow
\rcomod{A \tensor{B} A}$ establishes an equivalence of categories.
The converse holds if ${}_A(A \tensor{B} A)$ is flat.
\end{corollary}
\begin{proof}
Put $\Sigma = A$ in Theorem \ref{Descent}.
\end{proof}

\section{The structure of cosemisimple
corings}\label{cosemisimple}

Basic properties of cosemisimple corings have been studied in
\cite{ElKaoutit/Gomez/Lobillo:2001pr} and
\cite{Gomez/Louly:2001unp}. Perhaps, from the coring point of
view, the most fundamental examples of cosemisimple corings are
Sweedler's canonical corings $D \tensor{E} D$ for $E \subseteq D$
an extension of division rings. This section contains a full
description, in terms of finitely generated and projective right
$A$--modules and division subrings of their endomorphism rings, of
all cosemisimple $A$--corings for each fixed ring $A$.

 A coring is said to be \emph{cosemisimple} if it satisfies the
equivalent conditions in the following theorem.

\begin{theorem}\cite[Theorem 3.1]{ElKaoutit/Gomez/Lobillo:2001pr}\label{semi1}
Let $\coring{C}$ be an $A$--coring. The following statements are
equivalent:
\begin{enumerate}[(i)]
\item every left $\coring{C}$--comodule is semisimple and $\lcomod{\coring{C}}$ is abelian;
\item every right $\coring{C}$--comodule is semisimple and $\rcomod{\coring{C}}$ is is abelian;
\item $\coring{C}$ is semisimple as a left $\coring{C}$--comodule and $\coring{C}_A$ is flat;
\item $\coring{C}$ is semisimple as a right $\coring{C}$--comodule and ${}_A\coring{C}$ is flat;
\item $\coring{C}$ is semisimple as a right $\coring{C}^*$--module
and $\coring{C}_A$ is projective;
\item $\coring{C}$ is semisimple as a left
${}^*\coring{C}$--module and ${}_A\coring{C}$ is projective.
\end{enumerate}
\end{theorem}

This notion obviously generalizes cosemisimple coalgebras. On the
other hand, a ring $A$ is semisimple if and only if, considered as
$A$--coring, $A$ is cosemisimple. In fact, cosemisimple corings
were originally called semisimple corings in
\cite{ElKaoutit/Gomez/Lobillo:2001pr}, but it seems better, from
the point of view of the theory of Hopf algebras, to follow the
coalgebraic terminology.

Every cosemisimple $A$--coring $\coring{C}$ admits a unique
decomposition as a direct sum of simple cosemisimple
$A$--subcorings, where a coring is said to be \emph{simple} if it
has not non trivial subbicomodules \cite[Theorem
3.9]{ElKaoutit/Gomez/Lobillo:2001pr}. Several characterizations of
simple cosemisimple corings were given in \cite[Theorem
3.7]{ElKaoutit/Gomez/Lobillo:2001pr}. The structure of these
simple cosemisimple summands can be now deduced from our previous
results.

\begin{proposition}\label{simple1}
Let $\coring{C}$ be a simple cosemisimple $A$--coring and
$\Sigma_{\coring{C}}$ a finitely generated nonzero right
$\coring{C}$--comodule. Let $T = \rend{\coring{C}}{\Sigma}$ be the
simple artinian ring of endomorphisms of $\Sigma$. Then $\Sigma_A$
is finitely generated and projective and $can :
\rcomatrix{T}{\Sigma} \rightarrow \coring{C}$ is an isomorphism of
$A$--corings. Conversely, every comatrix $A$--coring
$\rcomatrix{B}{\Sigma}$, where $\Sigma_A$ is finitely generated
and projective and $B$ is simple artinian, becomes a simple
cosemisimple $A$--coring, and the number of simples in a complete
decomposition of $\Sigma_{\rcomatrix{B}{\Sigma}}$ coincides with
the number of simples in a complete decomposition of $B_B$.
\end{proposition}
\begin{proof}
The first statement is a consequence of Theorem \ref{generador}
because $\Sigma_{\coring{C}}$ is a finitely generated projective
generator of $\rcomod{\coring{C}}$. The second statement is a
consequence of Theorem \ref{Descent}.
\end{proof}

The following is our structure theorem for simple cosemisimple
corings.

\begin{theorem}\label{structuresimple}
An $A$--coring $\coring{C}$ is a simple cosemisimple coring if and
only if there is a finitely generated and projective right
$A$--module $\Sigma$ and a division subring $D \subseteq
\rend{A}{\Sigma}$ such that $\coring{C} \cong
\rcomatrix{D}{\Sigma}$ as $A$--corings. Moreover, if $\Xi$ is
another finitely generated and projective right $A$--module and $E
\subseteq \rend{A}{\Xi}$ is a division subring, then $\coring{C}
\cong \rcomatrix{E}{\Xi}$ as $A$--corings if and only if there is
an isomorphism of right $A$--modules $g : \Sigma \rightarrow \Xi$
such that $g D g^{-1} = E$.
\end{theorem}
\begin{proof}
That the simple cosemisimple $A$--corings are, up to isomorphism,
the comatrix corings $\rcomatrix{D}{\Sigma}$ for $\Sigma_A$
finitely generated and projective and $D \subseteq
\rend{A}{\Sigma}$ a division subring, is a consequence of
Proposition \ref{simple1}. Now, assume an isomorphism of
$A$--corings $f:\rcomatrix{D}{\Sigma} \rightarrow
\rcomatrix{E}{\Xi}$ as stated, and let $(-)_f:
\rcomod{\rcomatrix{D}{\Sigma}} \rightarrow
\rcomod{\rcomatrix{E}{T}}$ the associated induction functor
\cite[5.2]{Gomez:2002}, which is an isomorphism of categories. By
Proposition \ref{simple1}, $\Sigma$ is a simple right
$\rcomatrix{D}{\Sigma}$--comodule, which implies that $\Sigma_f$
is a simple right $\rcomatrix{E}{\Xi}$--comodule. Since $\Xi$ is,
up to isomorphism, the only simple comodule over
$\rcomatrix{E}{\Xi}$, there is an isomorphism of comodules $g :
\Sigma_f \rightarrow \Xi$, which, at the level of right
$A$--modules, gives an isomorphism $g : \Sigma_A \rightarrow
\Xi_A$. Clearly, $\rend{\rcomatrix{E}{\Xi}}{(\Sigma_f)} =
\rend{\rcomatrix{D}{\Sigma}}{\Sigma}$. By Theorem \ref{Descent},
$D = \rend{\rcomatrix{D}{\Sigma}}{\Sigma}$ and $E =
\rend{\rcomatrix{E}{\Xi}}{\Xi}$. Thus, if $d \in D$, then $g d
g^{-1}$ is an endomorphism of the comodule $\Xi$, that is, $g d
g^{-1} \in E$. We have obtained that $g D g^{-1} \subseteq E$. The
other inclusion is also easily obtained. Conversely, assume that
$g : \Sigma \rightarrow \Xi$ is an isomorphism of right
$A$--modules such that $g D g^{-1} = E$, and consider the
following $A$--bilinear map
\[
\psi: \Sigma^* \times \Sigma \rightarrow \rcomatrix{E}{\Xi},\qquad
((\varphi, u) \mapsto \varphi g^{-1} \tensor{E} g(u)).
\]
Let $d \in D$ and $e=gdg^{-1} \in E$, so $\psi(\varphi d,u)=
(\varphi d g^{-1} )\tensor{E} g(u) = \varphi  g^{-1}e \tensor{E}
g(u)= \varphi  g^{-1} \tensor{E} eg(u) = \varphi g^{-1} \tensor{E}
g(du) = \psi(\varphi,du)$. Hence $\psi$ extended to an
$A$--bilinear map $\psi: \rcomatrix{D}{\Sigma} \rightarrow
\rcomatrix{E}{\Xi}$. Given any right dual basis $\{e_i^*,e_i\}$
for $\Sigma_A$, it is easy to see that $\{e_i^* \circ g^{-1},
g(e_i)\}$ is a right dual basis for $\Xi_A$. Moreover, $\psi(e^*_i
\tensor{D}e_j) = (e^*_i \circ g^{-1}) \tensor{E} g(e_j)$, for all
$i,\,j$, thus $\psi$ is an isomorphism of $A$--bimodules. A direct
computation, using these bases, proves that $\psi$ is a morphism
of $A$--corings.
\end{proof}

Recall from \cite[Theorem 3.7]{ElKaoutit/Gomez/Lobillo:2001pr}
that an $A$--coring $\coring{C}$ is cosemisimple if and only if it
decomposes uniquely as $\coring{C} = \oplus_{\lambda \in \Lambda}
\coring{C}_{\lambda}$, where $\coring{C}_{\lambda}$ are simple
cosemisimple $A$--corings.

\begin{theorem}[Structure of cosemisimple corings]\label{structuresemisimple}
Let $\coring{C}$ be an $A$--coring. The following statements are
equivalents
\begin{enumerate}[(i)]
\item $\coring{C}$ is a cosemisimple $A$--coring;
\item there exists a family $\Lambda$ of finitely generated
projective right $A$--modules, and a division subring $D_{\Sigma}
\subseteq \mathrm{End}(\Sigma_A)$, for every $\Sigma \in \Lambda$,
such that $\coring{C} \cong \oplus_{\Sigma \in \Lambda} \left(
\Sigma^* \tensor{D_{\Sigma}} \Sigma \right)$ as $A$--corings.
\end{enumerate}
Furthermore, if $\coring{C}$ satisfies one of these conditions,
then the family $\Lambda$ is a set of representatives of all
simple right $\coring{C}$--comodules, and the decomposition is
unique in the following sense: given any other family $\Lambda'$
satisfying $(ii)$; then there exists a bijective map $\phi:
\Lambda \rightarrow \Lambda'$ and an isomorphism of right
$A$--modules $g_{\Sigma} : \Sigma_A \rightarrow \phi(\Sigma_A)$
for every $\Sigma \in \Lambda$ such that
$D_{\phi(\Sigma)}=g_{\Sigma} (D_{\Sigma}) g_{\Sigma}^{-1}$.
\end{theorem}
\begin{proof}
The equivalence $(i) \Leftrightarrow (ii)$ follows from
\cite[Theorem 3.7]{ElKaoutit/Gomez/Lobillo:2001pr} and Theorem
\ref{structuresimple}. Let us check the uniqueness of the
decomposition. Let $\Lambda$, $\Lambda'$ to be as in $(ii)$ with
the associated isomorphisms of $A$--corings
\[
\chi:\coring{C} \rightarrow \oplus_{\Sigma \in \Lambda} \left(
\Sigma^* \tensor{D_{\Sigma}} \Sigma \right), \qquad \zeta:
\coring{C} \rightarrow \oplus_{\Sigma' \in \Lambda'} \left(
\Sigma'{}^* \tensor{D_{\Sigma'}} \Sigma' \right).
\]
Therefore, $\coring{C} = \oplus_{\Sigma \in \Lambda}
\coring{C}_{\Sigma} = \oplus_{\Sigma' \in
\Lambda'}\coring{C}_{\Sigma'}$, where $\coring{C}_{\Sigma}=
\chi^{-1}(\rcomatrix{D_{\Sigma}}{\Sigma})$, $\coring{C}_{\Sigma'}=
\zeta^{-1}(\Sigma'{}^* \tensor{D_{\Sigma'}} \Sigma')$ are simple
cosemisimple $A$--corings for every $\Sigma \in \Lambda$ and
$\Sigma' \in \Lambda'$. The uniqueness of the decomposition given
in \cite[Theorem 3.7]{ElKaoutit/Gomez/Lobillo:2001pr} gives now
that there exists a bijection $\phi: \Lambda \rightarrow \Lambda'$
such that $\coring{C}_{\Sigma} = \coring{C}_{\phi(\Sigma)}$ for
every $\Sigma \in \Lambda$. That is,
$\rcomatrix{D_{\Sigma}}{\Sigma} \cong
\rcomatrix{D_{\phi(\Sigma)}}{\phi(\Sigma)}$, as $A$--corings, for
every $\Sigma \in \Lambda$. This implies in view of Theorem
\ref{structuresimple} that there exists an $A$--linear isomorphism
$g_{\Sigma}: \Sigma \rightarrow \phi(\Sigma)$ such that
$D_{\phi(\Sigma)} = g_{\Sigma} (D_{\Sigma}) g_{\Sigma}^{-1}$.
\end{proof}

\begin{remark}
The structure of cosemisimple coalgebras over a field $k$ is very
well-known: form direct sums of dual coalgebras of finite
dimensional simple algebras over $k$, which are matrices over
division $k$--algebras. Comatrix corings allow to built these
simple blocks directly, as they are coalgebras of the form
$\rcomatrix{D_{\Sigma}}{\Sigma}$, where the $D_{\Sigma}$'s are
division subalgebras of endomorphism algebras of finite
dimensional $k$--vector spaces $\Sigma$. Observe that this
description applies for coalgebras over arbitrary commutative
rings $k$, if we take finitely generated and projective
$k$--modules instead of finite dimensional vector spaces.
\end{remark}

\section{Coendomorphism corings}\label{coend}

We will see that comatrix corings are special instances of
coendomorphism corings. This gives an alternative approach,
although less elementary, for introducing comatrix corings and the
canonical map.

Let $\coring{C}$ and $\coring{D}$ be an $A$--coring and
$B$--coring, respectively. Let $N$ be an $A-B$--bimodule with a
right $\coring{D}$--coaction map $\rho_N : N \rightarrow N
\tensor{B} \coring{D}$ which is left $A$--linear. Assume that
$N_{\coring{D}}$ is quasi-finite, that is, the functor $-
\tensor{A} N : \rmod{A} \rightarrow \rcomod{\coring{D}}$ has a
left adjoint $F : \rcomod{\coring{D}} \rightarrow \rmod{A}$ (see
\cite[Section 3] {Gomez:2002}). This functor is called the
\emph{cohom} functor by analogy with the case of coalgebras over
fields (see \cite{Takeuchi:1977}); notation $F =
\cohom{\coring{D}}{N}{-}$. Let $\eta_{-,-} : \hom{\coring{D}}{-}{-
\tensor{A} N } \rightarrow \hom{A}{F(-)}{-}$ denote the natural
isomorphism of the adjunction, and $\theta :
1_{\rcomod{\coring{D}}} \rightarrow F(-) \tensor{A} N$ the unity
of the adjunction. The canonical map $A_A \rightarrow
\hom{\coring{D}}{N}{N} \rightarrow \hom{A}{F(N)}{F(N)}$ gives a
structure of left $A$--module on $F(N)$ such that $F(N)$ becomes
an $A$--bimodule. Define a comultiplication $\Delta : F(N)
\rightarrow F(N) \tensor{A} F(N)$ by $\Delta =
\eta_{N,F(N)\tensor{A}F(N)}\left((F(N) \tensor{A} \theta_N)\circ
\theta_N\right)$, that is, $\Delta$ is determined by the condition
$(F(N) \tensor{A} \theta_N)\circ \theta_N = (\Delta \tensor{A}
N)\circ \theta_N$; clearly $\Delta$ is $A$--bilinear. An analogous
proof to that of \cite[Proposition III.3.1]{AlTakhman:1999} shows
that $\Delta$ is coassociative. Moreover, $F(M)$ becomes an
$A$--coring with the counit given by $\epsilon =
\eta_{N,A}(\iota)$, where $\iota : N \rightarrow A \tensor{A} N$
is the canonical isomorphism. This $A$--coring will be denoted by
$e_{\coring{D}}(N)$; we refer to it as the \emph{coendomorphism}
$A$--coring associated to $N_{\coring{D}}$.

\begin{example}\label{quasi}
Let $\coring{C}$ be an $A$--coring, and $(\Sigma,\rho_{\Sigma})$ a
right $\coring{C}$--comodule such that $\Sigma_{A}$ is finitely
generated and projective with finite right dual basis
$\{e_i^*,e_i\}$; denote by $T=\rend{\coring{C}}{\Sigma}$ the
endomorphism ring of $\Sigma_{\coring{C}}$. It is easy to check
that the $A$--linear map
\[
\lambda_{\Sigma^*}: \Sigma^* \rightarrow \coring{C} \tensor{A}
\Sigma^*, \quad \left(\varphi \mapsto \sum_i((\varphi\tensor{A}
\coring{C})  \rho_{\Sigma})(e_i) \tensor{A} e^*_i \right),
\]
endows $\Sigma^*$ with a structure of left $\coring{C}$--comodule.
Moreover, $(-)^*: T \rightarrow \lend{\coring{C}}{\Sigma^*}$
sending $f \mapsto f^*$ establishes a ring isomorphism. Hence
$\Sigma^*$ is a $\coring{C}-T$--bicomodule.  There is a canonical
adjunction $-\tensor{T} \Sigma \dashv -\tensor{A} \Sigma^*$; the
associated natural isomorphism is
\[
\xymatrix@R=0pt@C=50pt{\eta_{Y_T,X_A}:
\hom{T}{Y}{X\tensor{A}\Sigma^*} \ar@{->}[r] &
\hom{A}{Y\tensor{T}\Sigma}{X} \\ f(y)=\sum_ix_i\tensor{A}e_i^*
\ar@{|->}[r] & (y\tensor{T}u \mapsto \sum_ix_ie_i^*(u)) \\ (y
\mapsto \sum_ig(y\tensor{T}e_i)\tensor{A} e_i^*) & g \ar@{|->}[l],
}
\]
the unit and counit are given by
\[
\begin{array}{ll}
\xymatrix@R=0pt{\theta_{Y_T}:Y \ar@{->}[r]& Y\tensor{T} \Sigma
\tensor{A}\Sigma^* \\  y \ar@{|->}[r] & \sum_iy\tensor{B}e_i
\tensor{A} e_i^*, } & \xymatrix@R=0pt{\chi_{X_A} : X \tensor{A}
\Sigma^* \tensor{T} \Sigma \ar@{->}[r] & X \\ x \tensor{A}e_i^*
\tensor{T} e_j \ar@{|->}[r]& xe_i^*(e_j). }
\end{array}
\]
Therefore, $\Sigma^*_T$ is quasi-finite (see \cite[Example
3.4]{Gomez:2002}). In this way, the coendomorphism $A$--coring
associated to $\Sigma^*_T$ is then $e_T(\Sigma^*)=\Sigma^*
\tensor{T} \Sigma$ with the following comultiplication and counit
\begin{multline*}
\Delta :e_T(\Sigma^*) \rightarrow e_{T}(\Sigma^*) \tensor{A}
e_T(\Sigma^*),\quad  (\varphi\tensor{T} u \mapsto \sum_{i,j}
\varphi \tensor{T} e_i \tensor{A} e^*_i \tensor{T} e_j e_j^*(u) \\
= \sum_i \varphi \tensor{T} e_i \tensor{A} e_i^* \tensor{T} u),
\end{multline*}
$\epsilon: \rcomatrix{T}{\Sigma} \rightarrow A$ sending $\varphi
\tensor{T} u \mapsto \sum_i\varphi(e_i)e_i^*(u) =\varphi(u)$. We
have shown that $e_T(\Sigma^*)$ is just the comatrix $A$--coring
of ${}_T\Sigma_A$.
\end{example}

Now, assume that $N$ is a $\coring{C}-\coring{D}$--bicomodule and
that ${}_B\coring{D}$ is a flat module. By \cite[Proposition
3.3]{Gomez:2002}, $F = \cohom{\coring{D}}{N}{-}$ factorizes
throughout the category $\rcomod{\coring{C}}$, and
$\cohom{\coring{D}}{N}{-} : \rcomod{\coring{D}} \rightarrow
\rcomod{\coring{C}}$ becomes a left adjoint to the cotensor
product functor $- \cotensor{\coring{C}} N : \rcomod{\coring{C}}
\rightarrow \rcomod{\coring{D}}$ with unity $\theta :
1_{\rcomod{\coring{D}}} \rightarrow F(-) \cotensor{\coring{C}} N$
and counity $\chi : F(- \cotensor{\coring{C}} N) \rightarrow
1_{\rcomod{\coring{C}}}$.

\begin{proposition}\label{cohom}
Let $N$ be a $\coring{C}-\coring{D}$--bicomodule that is
quasi-finite as a right $\coring{D}$--comodule, and assume that
${}_B\coring{D}$ is flat. The map $f : e_{\coring{D}}(N)
\rightarrow \coring{C}$ defined by $f = \chi_{\coring{C}} \circ
\cohom{\coring{D}}{N}{\lambda_N}$, where $\lambda_N : N
\rightarrow \coring{C} \tensor{A} N$ is the left comodule
structure map, is a homomorphism of $A$--corings.
\end{proposition}
\begin{proof}
The following diagram
\begin{equation}\label{diagrama1}
\xymatrix{ \coring{C} \cotensor{\coring{C}} N \ar@{=}[drr]
\ar^{\theta_{\coring{C} \cotensor{\coring{C}}N}}[rr] & &
F(\coring{C} \cotensor{\coring{C}} N) \cotensor{\coring{C}} N
\ar@{^{(}->}[r] \ar^{\chi_{\coring{C}} \cotensor{\coring{C}}
N}[d]& F(\coring{C} \cotensor{\coring{C}} N) \tensor{A} N
\ar^{\chi_{\coring{C}} \tensfun{A} N}[d]\\ & & \coring{C}
\cotensor{\coring{C}} N \ar@{^{(}->}[r] & \coring{C} \tensor{A}
N},
\end{equation}
where the hooked arrows are canonical monomorphisms, is
commutative, because the left triangle commutes by the properties
of the unity and counity of the adjunction. We have then that the
diagram
\begin{equation}
\xymatrix{F(N) \tensor{A} N \ar^{F(\lambda_N) \tensfun{A} N}[rr] &
& F(\coring{C} \cotensor{\coring{C}} N) \tensor{A} N
\ar^{\chi_{\coring{C}} \tensfun{A} N}[rr] & & \coring{C}
\tensor{A} N \\ N \ar^{\theta_N}[u] \ar^{\lambda_N}[rr] & &
\coring{C} \cotensor{\coring{C}} N \ar^{\theta_{\coring{C}
\cotensor{\coring{C}} N}}[u] \ar@{^{(}->}[urr] }
\end{equation}
is commutative. This means that $(f \tensor{A} N)\circ \theta_N =
\lambda_N$. To prove that $f$ is a homomorphism of $A$--corings,
we need to show that the following diagram is commutative:
\begin{equation}\label{diagrama2}
\xymatrix@C=60pt{F(N) \ar^{\Delta_{F(N)}}[r] \ar_{f}[d] & F(N)
\tensor{A} F(N) \ar^{f \tensfun{A} f}[d]\\ \coring{C}
\ar^{\Delta_{\coring{C}}}[r] & \coring{C} \tensor{A} \coring{C}},
\end{equation}
which is equivalent, by the adjunction isomorphism, to prove that
the diagram
\begin{equation}\label{diagrama3}
\xymatrix@C=45pt{N \ar^{\theta_N}[dr] & & & \\
 & F(N) \tensor{A} N \ar_{f \tensfun{A} N}[d]
 \ar^{F(N) \tensfun{A} \theta_N}[rr] & & F(N) \tensor{A} F(N) \tensor{A} N
 \ar^{f \tensfun{A} f \tensfun{A} N}[d]\\
  & \coring{C} \tensor{A} N \ar^{\Delta_{\coring{C}} \tensfun{A} N}[rr] & &
 \coring{C} \tensor{A} \coring{C} \tensor{A} N}
\end{equation}
is commutative. So
\[
(\Delta_{\coring{C}} \tensor{A} N)\circ (f \tensor{A} N)\circ
\theta_N = (\Delta_{\coring{C}} \tensor{A} N) \circ \lambda_N =
(\coring{C} \tensor{A} \lambda_N)\circ \lambda_N;
\]
\begin{multline*}
(f \tensor{A} f \tensor{A} N)\circ (F(N) \tensor{A} \theta_N)\circ
\theta_N = (f \tensor{A} ((f \tensor{A} N) \circ \theta_N))\circ
\theta_N \\ = (f \tensor{A} \lambda_N)\circ \theta_N = (\coring{C}
\tensor{A} \lambda_N)\circ (f \tensor{A} N)\circ \theta_N  =
(\coring{C} \tensor{A} \lambda_N)\circ \lambda_N
\end{multline*}
Therefore, \eqref{diagrama3} is commutative, and so is
\eqref{diagrama2}. Finally, we have to show that
$\epsilon_{\coring{C}}\circ f = \epsilon_{F(N)}$, which is
equivalent to show that $(\epsilon_{\coring{C}} \tensor{A} N)\circ
(f \tensor{A} N)\circ \theta_N = (\epsilon_{F(N)} \tensor{A}
N)\circ \theta_N$. This is clear from the commutative diagram
\[
\xymatrix@R=30pt@C=50pt{N \ar^{\theta_N}[dr] \ar^{\lambda_N}[drr]
\ar_{\iota}[ddr] & & \\
 & F(N) \tensor{A} N \ar_{f \tensfun{A} N}[r] \ar^{\epsilon_{F(N)}
  \tensfun{A} N}[d]& \coring{C}
 \tensor{A} N \ar^{\epsilon_{\coring{C}} \tensfun{A} N}[d] \\
  & A \tensor{A} N \ar@{=}[r] & A \tensor{A} N}
\]
\end{proof}

\begin{example}\label{pi=f}
Let $\coring{C}$, $\Sigma$ and $T$ as in Example \ref{quasi}.
Then, by Proposition \ref{cohom}, there exists an $A$--coring
morphism $f = \chi_{\coring{C}} \circ (\lambda_{\Sigma^*}
\tensor{T} \Sigma) :e_T(\Sigma^*) \rightarrow \coring{C}$. Using
the notation of Example \ref{quasi}, we have $f(\varphi\tensor{T}
u)= \sum_i \chi_{\coring{C}}( ((\varphi \tensor{A}
\coring{C})\circ \rho_{\Sigma}(e_i)) \tensor{A} e_i^* \tensor{T} u
)= \sum_i (\varphi \tensor{A}
\coring{C})\rho_{\Sigma}(e_i)e_i^*(u) = (\varphi \tensor{A}
\coring{C}) \circ \rho_{\Sigma}(u)$. Therefore, $f=can$, the
morphism of $A$--corings defined in Proposition \ref{pi}.
\end{example}

 \providecommand{\bysame}{\leavevmode\hbox
to3em{\hrulefill}\thinspace}

\end{document}